\documentclass[tikz,border=5pt]{article}
\usepackage{etoolbox}
\usepackage{amssymb}
\usepackage{amsmath,amsthm}
\usepackage{amsmath} 
\DeclareMathOperator{\codeg}{codeg}
\providecommand{\dotcup}{\mathbin{\dot{\cup}}}
\usepackage{mathtools}
\usepackage{bm}
\usepackage{extarrows}
\usepackage{indentfirst}
\usepackage[english]{babel}
\usepackage{tikz}
\usepackage{float}
\usepackage{lmodern} 
\usetikzlibrary{calc}
\usetikzlibrary{positioning,shapes,arrows.meta, shapes.geometric,fit}
\newcommand{\ex}{\mathrm{ex}}
\usepackage{authblk}
\providecommand{\Ik}[1]{\overline{K_{#1}}}

\providecommand{\Z}{\mathbb{Z}}
\DeclareMathOperator{\dist}{dist}

\newcommand{\hk}{h_k}

\providecommand{\gk}{g_k}
\theoremstyle{thmstyleone}%
\newtheorem{theorem}{Theorem}[section]
\newtheorem{proposition}[theorem]{Proposition}%

\theoremstyle{thmstyletwo}%

\theoremstyle{thmstylethree}%
\newtheorem{lemma}{Lemma}

\newtheorem{thm}{Theorem}[section]

\newtheorem{conj}[thm]{Conjecture}

\newtheorem{claim}{Claim}[section]

\newtheorem{corollary}{Corollary}

\title{The generalized Tur\'an number for $K_3$ in graphs without suspensions of a path 
}
\author[1]{Doudou Hei}
\author[2,3]{Xinmin Hou}

\affil[1]{School of Mathematical Sciences, Suzhou University of Science and Technology, Suzhou, Jiangsu 215009, China}
\affil[2]{School of Mathematical Sciences, University of Science and Technology of China,  Hefei, Anhui 230026, China}
\affil[3]{Hefei National Laboratory, University of Science and Technology of China,  Hefei, Anhui 230088, China}

\date{}

\begin{document}

\maketitle
\begin{abstract}
  Given graphs $H$ and $F$, the generalized Tur\'an number $\ex(n, H, F)$ is defined as 
  the maximum number of copies of $H$ in an $n$-vertex graph that contains no copy of $F$.
  The suspension $\widehat{F}$ of a graph $F$ is 
  obtained  by adding a new vertex that is adjacent to every vertex of $F$.
  Mubayi and Mukherjee (2023, DM) conjectured that $\ex(n, K_3, \widehat{P_k})=\left\lfloor \frac{k-2}{2}\right\rfloor \cdot \frac{n^2}{8}+o(n^2)$, where $P_k$ is a
  path on $k\ge 4$ vertices, and they verified this conjecture for $k=4,5,6$. 
  Recently, this conjecture was independently proved correct by Dou and Peng and by Wang.
  In this paper, using the stability method, we determine the exact value of $\ex(n, K_3, \widehat{P_k})$, where $k\geq 6$ for sufficiently large $n$.
\end{abstract}

\section{Introduction}\label{sec1}
All graphs considered in this paper are finite, simple, and
undirected. A graph is called $F$-free if it contains no copy
of $F$ as a subgraph, not necessarily induced. For a graph $G$, write $v(G)$ and $e(G)$ for its numbers of
vertices and edges.
For graphs $H$ and $F$, the generalized Tur\'an number is
defined by
$$\ex(n, H , F)=\max\{N(H, G): G \text{ is } F\text{-free}, \vert V(G)\vert=n\}.$$
where $N(H,G)$ is the number of unlabelled copies of $H$ in $G$.
We write $t(G)=N(K_3,G)$ for the number of triangles in $G$.
Let $P_l$, $C_l$
 $K_l$, and $\overline{K_l}$ denote a path, cycle, complete and empty graph on $l$ vertices, respectively. For vertex-disjoint graphs $H$ and $J$, their join $H+J$ is
obtained from their disjoint union by adding all edges between
$V(H)$ and $V(J)$.

The suspension  of a graph $F$, denoted by $\widehat{F}$, is the graph obtained by adding an additional vertex to $F$ and connecting it to every vertex of $F$.
Mubayi and Mukherjee \cite{7311} studied $\ex(n, K_3, \widehat{F})$ for various different bipartite graphs $F$. In particular, they investigated $\ex(n, K_3, \widehat{P}_k)$ and gave the following bound and open problem.
\begin{proposition}[Mubayi and Mukherjee \cite{7311}]\label{lemma 1}
  Let $n\geq k\geq 4$. Then 
  $$\left\lfloor \frac{k-2}{2}\right\rfloor \cdot \frac{n^2}{8}\leq \ex(n,K_3, \widehat{P_k})\leq \frac{k-2}{12}\cdot n^2+\frac{(k-2)^2}{12}\cdot n, $$
  where the lower bound holds when $n$ is a multiple of $4\left\lfloor \frac{k-2}{2}\right\rfloor$.
\end{proposition}

Mubayi and Mukherjee \cite{7311} believe that the lower bound above is asymptotically
tight for all fixed $k\geq 4$ and propose the following conjecture.\\

\begin{conj}[Mubayi and Mukherjee \cite{7311}]\label{Conj:MM23}
Let $n\geq k\geq 4$. Then 
  $$\ex(n,K_3, \widehat{P_k})=\left\lfloor \frac{k-2}{2}\right\rfloor \cdot \frac{n^2}{8}+o(n^2).$$
\end{conj}

They proved the conjecture for the 
first three cases: $k=4, 5$, and 6. 
\begin{theorem}[Mubayi and Mukherjee \cite{7311}]\label{thm 1}
  For $k=4, 5$ and $6$, 
  $$\ex(n,K_3, \widehat{P_k})=\left\lfloor \frac{k-2}{2}\right\rfloor \cdot \frac{n^2}{8}+o(n^2).$$
  When $k=4$ or $k=6$, the error term can be improved to $O(n)$.
\end{theorem}

Recently, Dou and Peng\cite{11191} determined the asymptotic value of $\ex\left(n, K_3, \widehat{P_k}\right)$ for every $k \geq 4$ using a unified method and thus prove Conjecture \ref{Conj:MM23}.
In addition, Mukherjee\cite{874} determined the exact value of $\operatorname{ex}\left(n, K_3, \widehat{P}_4\right)$ for all $n$, and a similar result was proved by Gerbner \cite{873} for large $n$. Hei, Hou, and Ma\cite{11192}  determined the exact value of $\operatorname{ex}\left(n, K_3, \widehat{P}_5\right)$ for sufficiently large $n$.

Here we consider fixed $k\ge 6$ and develop
a structural reduction that also identifies the linear
term and a finite optimization for the remaining bounded
correction. As a result, we determine the exact value of $\ex(n, K_3, \widehat{P_k})$ for sufficiently large $n$.
Throughout the paper, let
\begin{equation}
r=\left\lfloor\frac{k-2}{2}\right\rfloor,\qquad
q=r+1.
\end{equation}
Define
\begin{equation}\label{eq:g-definition-intro}
g_k(n)
=
\max_H\bigl\{(n-v(H))e(H)+t(H)\bigr\},
\end{equation}
where the maximum is over all $P_k$-free graphs $H$ satisfying
$v(H)\le n$ and $\Delta(H)\leq r$.

\begin{theorem}\label{thm:main-structure}
For every fixed integer $k\ge 6$, there exists $n_0=n_0(k)$
such that  for all $n\ge n_0$,
$$ex(n, K_3, \widehat{P_k})=g_k(n).$$
And every extremal graph is isomorphic to $
H+\overline K_{n-v(H)}$,
where $H$ is an $r$-regular $P_k$-free graph such that 
$g_k(n)$ attains its maximum value.
\end{theorem}

For the exact formula, we have the following result.
\begin{theorem}\label{main2}
For every fixed integer $k \geq 6$ and all sufficiently large $n$,

$$
ex(n, K_3, \widehat{P_k})=\frac{r}{8} n^2+\frac{r(r-1)}{12} n+c_k(n),
$$
where $c_k(n)$ is bounded and eventually periodic, with $2 q$ as a period.  
\end{theorem}

The rest of this article is arranged as follows. In Section 2, we present some preliminary results. The proof of Theorem
\ref{thm:main-structure} and Theorem \ref{main2} are given in Section 3 and 4, respectivly.

\section{Preliminary}\label{sec2}
In this section, we give preparatory lemmas for the proof of Theorem \ref{thm:main-structure}.
Let
$$
x_n=\frac{n}{2}+\frac{r-1}{6}, \quad U(n)=\frac{r}{8} n^2+\frac{r(r-1)}{12} n+\frac{r(r-1)^2}{72} .
$$
Fix a graph $G$. For $v\in V(G)$, let $N_G(v)$ be its
neighborhood, and set
\[
t_v(G)=e(G[N_G(v)]).
\]
Thus $t_v(G)$ is the number of triangles containing $v$.
For a vertex set $X\subseteq V(G)$, write
\[
N_X(v)=N_G(v)\cap X,\qquad d_X(v)=|N_X(v)|.
\]
For disjoint sets $X,Y$, let $e_G(X,Y)$ denote the number of
edges with one endpoint in each set.

We would repeatedly use the equivalence
\begin{equation}\label{eq:neighborhood-characterization}
G\text{ is }\widehat{P_k}\text{-free}
\quad\Longleftrightarrow\quad
G[N_G(v)]\text{ is }P_k\text{-free for every }v\in V(G).
\end{equation}

To begin with, we establish a lower bound on the number of triangles in extremal graphs.
\begin{lemma}\label{lem:join}
Suppose that $H$ is $P_k$-free and $\Delta(H)\le r$.
For every integer $b\ge 0$, the graph $H+\Ik{b}$ is $\widehat{P_k}$-free,
and
\[
t(H+\Ik{b})=b\,e(H)+t(H).
\]
\end{lemma}

\begin{proof}
The neighborhood of a vertex in $\Ik{b}$ induces $H$ and is
therefore $P_k$-free.
For $u\in V(H)$, its neighborhood consists of $N_H(u)$
and the independent set $V(\Ik{b})$.
Every path in this neighborhood contains at most
$|N_H(u)|$ vertices from $H$ and at most $|N_H(u)|+1$
vertices from $\Ik{b}$.
Its order is consequently at most
\[
2|N_H(u)|+1\le 2r+1<k.
\]
The first assertion follows from
\eqref{eq:neighborhood-characterization}.
Every triangle is either contained in $H$ or consists of
an edge of $H$ and a vertex of $\Ik{b}$, proving the formula.
\end{proof}
Taking $H=(a/q)K_q$, where $q\mid a$, gives
\[
t(H+\Ik{n-a})
=
\frac r2a(n-a)+\frac{r(r-1)}6a.
\]
Define
\begin{equation}\label{eq:h-definition}
\hk(n)
=
\max_{\substack{0\le a\le n\\q\mid a}}
\left\{
\frac r2a(n-a)+\frac{r(r-1)}6a
\right\}.
\end{equation}
Completing the square gives
\begin{equation}\label{eq:complete-square}
\frac r2a(n-a)+\frac{r(r-1)}6a
=
U(n)-\frac r2(a-x_n)^2.
\end{equation}
For all sufficiently large $n$, there is a multiple of $q$
in $[0,n]$ at distance at most $q/2$ from $x_n$.
Therefore
\begin{equation}\label{eq:h-envelope}
U(n)-\frac{rq^2}{8}\le \hk(n)\le U(n)
\end{equation}
for all sufficiently large $n$, and
\begin{equation}\label{eq:lower-bound}
\ex(n, K_3, \widehat{P_k})\ge \hk(n)
=
\frac r8n^2+\frac{r(r-1)}{12} n+O_k(1).
\end{equation}

The Tur\'an number for path is given by Erdős and Gallai~\cite{875}.
\begin{theorem}[Erdős and Gallai \cite{875}]\label{thm8121}
  $$\ex(n, P_k)\leq \frac{n}{k-1}\binom{k-1}{2}= \frac{k-2}{2}n. $$
\end{theorem}

Now we can obtain an inequality relating the number of triangles and the number of edges in an extremal graph. The important tools for proving
it are the following result from Dou and Peng\cite{11191} and  the triangle removal lemma\cite{ruzsa1978triple}.
\begin{lemma}[Dou and Peng\cite{11191}]\label{lem:DP}
There exists a constant $C=C(k)\ge 0$ such that the following
holds. Let $Q$ be a $P_k$-free graph with a partition
$V(Q)=X\cup Y$. If
\[
d_X(y)\le r
\quad\text{and}\quad
d_Y(y)\le r
\qquad\text{for every }y\in Y,
\]
then
\begin{equation}\label{eq:DP}
e_Q(X,Y)+2e(Q[Y])\le r|Y|+C|X|.
\end{equation}
\end{lemma}

It is easy to obtain a one-sided bounded-degree version from the above lemma.
\begin{lemma}\label{lem:one-sided}
There exists a constant $D=D(k)\ge 0$ such that, whenever
$Q$ is $P_k$-free with a partition $V(Q)=X\cup Y$ and
$\Delta(Q[Y])\le r$, then
\begin{equation}\label{eq:one-sided}
e_Q(X,Y)+2e(Q[Y])\le r|Y|+D|X|.
\end{equation}
\end{lemma}

\begin{proof}
Put
\[
Y^+=\{y\in Y:d_X(y)\ge r+1\},\qquad
Y^-=Y\setminus Y^+,\qquad
\alpha=\frac{k-2}{2}.
\]
The bipartite graph consisting of the edges between $X$
and $Y^+$ is $P_k$-free. Hence
\[
(r+1)|Y^+|
\le e_Q(X,Y^+)
\le \alpha(|X|+|Y^+|).
\]
Since $\alpha\in\{r,r+\tfrac12\}$,
\begin{equation}\label{eq:Yplus}
|Y^+|
\le \frac{\alpha}{r+1-\alpha}|X|
\le (2r+1)|X|,
\end{equation}
and then we have,
\[
e_Q(X,Y^+)\le \alpha(2r+2)|X|.
\]

Apply Lemma~\ref{lem:DP} to $Q[X\cup Y^-]$, we get 
$$e_Q(X,Y^-)+2e(Q[Y^-])\le r|Y^-|+C|X|.$$
The contribution to $2e(Q[Y])$ from edges meeting $Y^+$
is at most $2r|Y^+|$.
Consequently,
\[
e_Q(X,Y)+2e(Q[Y])
\le
r|Y^-|+
\bigl(C+\alpha(2r+2)+2r(2r+1)\bigr)|X|.
\]
Let
\[
D=C+\alpha(2r+2)+2r(2r+1),
\] and this proves the result.
\end{proof}

\begin{lemma}[The triangle Removal Lemma\cite{ruzsa1978triple}]
  For each $\epsilon>0$, there is $\delta=\delta(\epsilon)>0$ such that every graph on $n$ vertices with at most $\delta n^3$ copies of $K_3$ can be made $K_3$-free by removing at most $\epsilon n^2$ edges.
\end{lemma}

\begin{theorem}\label{thm:triangle-edge}
For every fixed $k\ge 6$ and every $\varepsilon>0$,
there exists $N=N(k,\varepsilon)$ such that every
$n$-vertex $\widehat{P_k}$-free graph $G$ with $n\ge N$ satisfies
\begin{equation}\label{eq:triangle-edge}
t(G)\le \frac r2e(G)+\varepsilon n^2.
\end{equation}
\end{theorem}

\begin{proof}
By \eqref{eq:neighborhood-characterization} and
Theorem \ref{thm8121},
\begin{equation}\label{eq:crude-triangles}
3t(G)
=
\sum_{v\in V(G)}e(G[N_G(v)])
\le
\frac{k-2}{2}\sum_{v\in V(G)}d_G(v)
=
(k-2)e(G).
\end{equation}
In particular,
\[
t(G)\le \frac{k-2}{6}n^2=o(n^3).
\]
Fix $\rho\le \frac{\varepsilon^2}{(2(k-2)+C)\varepsilon+4(2r+1)(k-2)^2}$, where $C$ is from Lemma \ref{lem:DP}. By the triangle removal lemma
\cite{ruzsa1978triple}, for sufficiently large $n$ there is a set
$R_0\subseteq E(G)$ such that
\[
|R_0|\le \rho n^2
\quad\text{and}\quad
G-R_0\text{ is triangle-free}.
\]

Let $W=\frac{(k-2)^2}{\varepsilon}$,  and define
\[
R_1=
\{uv\in E(G)\setminus R_0:\codeg_G(uv)>W\},
\qquad
R=R_0\cup R_1,
\qquad
R_3=E(G)\setminus R.
\]
Here $\codeg_G(uv)=|N_G(u)\cap N_G(v)|$.
By \eqref{eq:crude-triangles},
\begin{equation}\label{eq:R1-size}
|R_1|
\le
\frac1W\sum_{uv\in R_1}\codeg_G(uv)
\le \frac{3t(G)}W
\le \frac{k-2}{2W}n^2.
\end{equation}
Thus
\begin{equation}\label{eq:R-size}
|R|\le
\left(\rho+\frac{k-2}{2W}\right)n^2.
\end{equation}

For an edge set $S$, let
\[
N_S(v)=\{u:uv\in S\}.
\]
Let $\mathcal D$ be the family of triangles containing
at least two edges of $R$.
Every such triangle is counted at least once in
$\sum_v e(G[N_R(v)])$. Therefore
\begin{equation}\label{eq:two-R}
|\mathcal D|
\le
\sum_v e(G[N_R(v)])
\le
\frac{k-2}{2}\sum_v|N_R(v)|
=
(k-2)|R|.
\end{equation}

For an edge $vy\in R_3$, define the oriented quantity
\[
c_v(vy)=|N_{R_0}(v)\cap N_{R_3}(y)|.
\]
Call $vy$ bad if $c_v(vy)\ge r+1$ or $c_y(vy)\ge r+1$,
and let $\mathcal B$ be the set of bad edges.
Fix $v$, let
\[
X=N_{R_0}(v),\qquad
Y=N_{R_3}(v),\qquad
Y^+=\{y\in Y:c_v(vy)\ge r+1\}.
\]
The sets $X$ and $Y$ are disjoint.
Consider the bipartite subgraph between
$X$ and $Y^+$ in $G[N_G(v)]$, it is $P_k$-free, and we obtain
\[
(r+1)|Y^+|
\le e_G(X,Y^+)
\le \frac{k-2}{2}(|X|+|Y^+|).
\]
Then we have 
\[
|Y^+|\le (2r+1)|X|.
\]
Summing over $v$ gives
\begin{equation}\label{eq:bad-edges}
|\mathcal B|
\le 2(2r+1)|R_0|.
\end{equation}
Since every edge of $\mathcal B$ has codegree at most $W$,
the number of triangles containing a bad edge is at most
\begin{equation}\label{eq:bad-triangles}
2(2r+1)W|R_0|.
\end{equation}

Let $R_4=R_3\setminus\mathcal B$, and let
$\mathcal U$ be the family of triangles with one edge in
$R_0$ and two edges in $R_4$.
For each $v$, define
\[
X_v=N_{R_0}(v),\qquad Y_v=N_{R_4}(v).
\]
Let $Q_v$ have vertex set $X_v\dotcup Y_v$,
with the following edges:
\begin{itemize}
\item edges of $R_4$ between $X_v$ and $Y_v$;
\item edges of $R_0$ with both endpoints in $Y_v$.
\end{itemize}
The graph $Q_v$ is a subgraph of $G[N_G(v)]$ and is
therefore $P_k$-free.

For $y\in Y_v$, the edge $vy$ belongs to $R_4$.
Hence
\[
d_{Q_v,X_v}(y)\le c_v(vy)\le r,
\qquad
d_{Q_v,Y_v}(y)\le c_y(vy)\le r.
\]
Lemma~\ref{lem:DP} gives
\[
e_{Q_v}(X_v,Y_v)+2e(Q_v[Y_v])
\le r|Y_v|+C|X_v|.
\]
Upon summing over $v$, each triangle of $\mathcal U$
contributes once at each endpoint of its $R_0$-edge and
twice at its third vertex. Thus
\[
4|\mathcal U|
\le 2r|R_4|+2C|R_0|,
\]
and consequently
\begin{equation}\label{eq:U-bound}
|\mathcal U|
\le \frac r2|R_4|+\frac C2|R_0|
\le \frac r2e(G)+\frac C2|R_0|.
\end{equation}

Every triangle of $G$ contains an edge of $R_0$.
If it is not in $\mathcal D$, its other two edges lie in
$R_3$; if neither is bad, it belongs to $\mathcal U$.
Combining \eqref{eq:R-size}, \eqref{eq:two-R},
\eqref{eq:bad-triangles}, and \eqref{eq:U-bound}, we obtain
\[
t(G)
\le \frac r2e(G)
+
\left[
\frac{(k-2)^2}{2W}
+
\left(k-2+2(2r+1)W+\frac C2\right)\rho
\right]n^2\\
\le \frac r2e(G)+\varepsilon n^2.
\]
\end{proof}

From the above theorem, we immediately obtain a result of enumerating the edges of extremal configurations.

\begin{corollary}\label{11211}
 For every $\varepsilon >0$, there exists $n_0=n_0(\varepsilon)$ such that the following holds. If $G$ is a $\widehat{P_k}$-free  graph on $n>n_0$ vertices with $t(G)\geq \frac{rn^2}{8}-\varepsilon n^2$, 
  then $e(G)\geq \frac{n^2}{4}-\frac{4}{r}\varepsilon n^2.$
\end{corollary}

Next, we employ the following theorem to carry out the preliminary characterization of extremal graphs.
Let $T(n, k)$ denote the Tur\'an graph on $n$ vertices with $k$ partites.
\begin{theorem}[Erdős-Simonovits Stability Theorem \cite{821}]\label{thm89} Fix a graph $F$ with $\chi(F) \geq 3$.
For any $\varepsilon>0$, there exist $\delta=\delta(\varepsilon )>0$ and $n_1=n_1(\varepsilon )$ such that if $G$ is an $n>n_1$ vertices $F$-free graph with

$$
e(G)>\operatorname{ex}(n, F)-\delta n^2,
$$
then the edit distance between $G$ and $T(n, \chi(F)-1)$ is at most $\varepsilon n^2$. In other words, we can add and delete at most $\varepsilon n^2$ edges of $G$ to obtain $T(n, \chi(F)-1)$.
  
\end{theorem}

The extremal number of $\widehat{P_k}$ is determined by Zhu, Wang, Zhang, and Zhang \cite{11202}.
\begin{theorem}[Zhu, Wang, Zhang, and Zhang \cite{11202}]\label{11261}
When $n \geq 4(4 k)^6$,
$$\ex(n, \widehat{P_k})\leq\frac{n^2}{4}+\left\lfloor \frac{n+1}{4}\right\rfloor.$$
\end{theorem}

By Corollary \ref{11211}, Theorem \ref{thm89} and \ref{11261}, we have the following result.
\begin{corollary}\label{cor:stability}
Let $G$ be an $n$-vertex $F_k$-free graph satisfying
$t(G)\ge \hk(n)$.
Then, as $n\to\infty$,
there is a partition $V(G)=V_1\cup V_2$ such that
\begin{gather}
|V_1|,|V_2|=\frac n2+o(n),\label{eq:balanced}\\
e(G[V_1])+e(G[V_2])=o(n^2),\label{eq:internal-small}\\
|V_1||V_2|-e_G(V_1,V_2)=o(n^2).\label{eq:missing-small}
\end{gather}
\end{corollary}


\section{Proof of Theorem \ref{thm:main-structure}}

First, we employ a standard technique to refine the extremal graph.
Recall the definition of $\gk(n)$ in
\eqref{eq:g-definition-intro}. By Lemma~\ref{lem:join},
\begin{equation}\label{eq:h-g-ex}
\hk(n)\le \gk(n)\le \ex(n, K_3, \widehat{P_k}).
\end{equation}

\begin{proposition}\label{prop:g-envelope}
For every $n\ge 0$,
\[
\hk(n)\le \gk(n)\le U(n).
\]
Moreover,
\[
\gk(n)=\frac r8n^2+\frac{r(r-1)}{12} n+O_k(1).
\]
\end{proposition}

\begin{proof}
Let $H$ be admissible in \eqref{eq:g-definition-intro},
Since $\Delta(H)\le r$,
\[
e(H)\le \frac{rv(H)}{2},
\qquad
3t(H)=\sum_{v\in V(H)}e(H[N_H(v)])
\le v(H)\binom r2.
\]
It follows that
\begin{align*}
(n-v(H))e(H)+t(H)
&\le \frac{rv(H)}2(n-v(H))+\frac{r(r-1)}6v(H) \\
&= U(n)-\frac{r}{2}(v(H)-x_n)^2 \le U(n).
\end{align*}
The asymptotic formula follows from
\eqref{eq:h-envelope}.
\end{proof}

\begin{lemma}\label{lem:deletion}
There is a constant $C_0=C_0(k)>0$ with the following property.
If $G_n$ is an $n$-vertex $\widehat{P_k}$-free graph satisfying
$t(G_n)\ge \gk(n)$, then $G_n$ has an induced subgraph $G_m$
such that
\begin{align}
t(G_m)&\ge \gk(m)+n-m,\label{eq:deletion-surplus}\\
t_v(G_m)&\ge \frac r4m-C_0,
v\in V(G_m).\label{eq:minimum-triangle-degree}
\end{align}
Furthermore, $m\to\infty$ as $n\to\infty$.
\end{lemma}

\begin{proof}
By Proposition~\ref{prop:g-envelope}, 
we may assume that
\begin{equation}\label{eq:g-increment}
\gk(j)-\gk(j-1)\ge \frac r4j-C_0+1
\qquad(j\ge 1).
\end{equation}
Indeed, for large $j$ this follows by subtracting the upper
bound $U(j-1)$ from the lower bound
$U(j)-rq^2/8$.

Start with $G_n$.
Whenever the current graph $G_j$ contains a vertex $v$ with
\[
t_v(G_j)<\frac r4j-C_0,
\]
delete $v$ and let $G_{j-1}=G_j-v$.
If $t(G_j)\ge \gk(j)+n-j$, then
\[
\begin{aligned}
t(G_{j-1})
&>
\gk(j)+n-j-\left(\frac r4j-C_0\right)\\
&\ge
\gk(j-1)+n-(j-1),
\end{aligned}
\]
by \eqref{eq:g-increment}.
Thus the surplus inequality is preserved.
At termination, both
\eqref{eq:deletion-surplus} and
\eqref{eq:minimum-triangle-degree} hold.

Finally, $\gk(m)\ge 0$, so
\[
n-m\le t(G_m)\le \binom m3.
\]
Hence $m$ cannot remain bounded while $n\to\infty$.
\end{proof}
In the remainder of this section, $G=G_m$ is a graph obtained
from Lemma~\ref{lem:deletion}.
By Theorem \ref{thm8121},
\[
t_v(G)\le \frac{k-2}{2}d_G(v).
\]
Consequently,
\begin{equation}\label{eq:minimum-degree}
\delta(G)\ge \frac{r}{2(k-2)}m-O_k(1).
\end{equation}
Also,
\[
t(G)\ge \gk(m)\ge \hk(m),
\]
so Corollary \ref{cor:stability} applies.
By 
Choose a maximum cut $V(G)=V_1\dotcup V_2$.
It satisfies

\begin{gather}
d_{V_3-i}(v)\geq d_{V_i}, v\in V_i,\label{eq:max-cut}\\
|V_1|,|V_2|=\frac m2+o(m),\label{eq:balanced}\\
e(G[V_1])+e(G[V_2])=o(m^2),\label{eq:internal-small}\\
|V_1||V_2|-e_G(V_1,V_2)=o(m^2).\label{eq:missing-small}
\end{gather}


By \eqref{eq:minimum-degree} and \eqref{eq:max-cut}, for all
sufficiently large $m$,
\begin{equation}\label{eq:cross-minimum}
d_{V_{3-i}}(v)\ge \frac{r}{8(k-2)}m
\qquad(v\in V_i).
\end{equation}

Fix
\begin{equation}\label{eq:explicit-eta}
\eta=\eta(k)=\frac{1}{128(k-2)(r+1)^2},\quad c=\frac{r}{8(k-2)}.
\end{equation}
This choice satisfies
\begin{equation}\label{eq:eta-conditions}
(r+1)\eta<\frac{c}{2},\qquad
k\eta<\frac14,\qquad
r(r+1)\eta\le \frac1{32}.
\end{equation}

For $i\in\{1,2\}$, define
\[
E_i=
\{v\in V_i:|V_{3-i}\setminus N_G(v)|>\eta m\},
\qquad
T_i=V_i\setminus E_i.
\]
By \eqref{eq:missing-small},
\begin{equation}\label{eq:E-small}
|E_1|+|E_2|=o(m).
\end{equation}

\begin{claim}\label{clm:good-degree}
Every vertex $v\in V_i$
has at most $r$ neighbors in $T_i$.
In particular, $\Delta(G[T_i])\le r$.
\end{claim}

\begin{proof}
Suppose that $v$ has distinct neighbors
$u_1,\ldots,u_{r+1}\in T_i$.
By \eqref{eq:cross-minimum}, the common neighbors of
$v,u_1,\ldots,u_{r+1}$ in $V_{3-i}$ number at least
\[
cm-(r+1)\eta m.
\]
By \eqref{eq:eta-conditions}, this is at least $r+2$ for
sufficiently large $m$.
Thus $G[N_G(v)]$ contains $K_{r+1,r+2}$.
This complete bipartite graph contains a path on $2r+3$
vertices and hence contains $P_k$, a contradiction.
\end{proof}

\begin{claim}\label{clm:good-path-free}
The graph $G[T_i]$ is
$P_k$-free.
\end{claim}
\begin{proof}
Suppose that $G[T_i]$ contains a copy of $P_k$.
Its $k$ vertices have at least
\[
|V_{3-i}|-k\eta m
=
\frac m2-k\eta m+o(m)>0
\]
common neighbors in $V_{3-i}$.
Any such common neighbor, together with the path, gives
a copy of $\widehat{P_k}$, a contradiction.
\end{proof}

Let
\begin{equation}\label{eq:closure-constants}
L_k=\sum_{j=0}^{k-2}r^j,
\qquad
K_k=1+rL_k.
\end{equation}
Every connected component of $G[T_i]$ has at most $L_k$
vertices. To see this, root a breadth-first search tree in
such a component. Its depth is at most $k-2$, since a longer
root-to-vertex path would contain $P_k$, and its maximum
branching is at most $r$.

Let $B_i$ consist of $E_i$ together with all components of
$G[T_i]$ that meet $N_G(E_i)\cap T_i$.
Set
\[
A_i=V_i\setminus B_i,\qquad
a_i=|A_i|,\qquad b_i=|B_i|,\qquad b=b_1+b_2,\quad B=B_1\cup B_2
\]
Claim~\ref{clm:good-degree} gives
\[
|N_G(E_i)\cap T_i|\le r|E_i|.
\]
Therefore
\begin{equation}\label{eq:B-size}
|B_i|\le K_k|E_i|,
\qquad b=o(m).
\end{equation}
By construction,
\begin{equation}\label{eq:clean-properties}
e_G(A_i,B_i)=0,\qquad
\Delta(G[A_i])\le r,\qquad
G[A_i]\text{ is }P_k\text{-free}.
\end{equation}
Also,
\begin{equation}\label{eq:A-balanced}
a_i=\frac m2+o(m).
\end{equation}

Let  $\tau_B(G)$ be the number of
triangles of $G$ meeting $B$.

\begin{claim}\label{clm:exceptional-triangles}
\begin{equation}\label{eq:exceptional-bound}
\tau_B(G)
\le
\frac r2
\sum_{i=1}^2\sum_{u\in B_i}d_{A_{3-i}}(u)
+O_k(b^2).
\end{equation}
\end{claim}

\begin{proof}
For $u\in B_i$, put
\[
X_u=N_G(u)\cap B_i,\qquad
Y_u=N_G(u)\cap A_{3-i}.
\]
The graph $G[X_u\cup Y_u]$ is $P_k$-free, and $\Delta(G[Y_u])\le r$.
Lemma~\ref{lem:one-sided} gives
\[
e(G[Y_u])+\frac12e_G(X_u,Y_u)
\le \frac r2|Y_u|+\frac D2|X_u|.
\]

Since $e_G(A_i,B_i)=0$, every triangle meeting $B$ but
not contained in $B$ either has one vertex in $B_i$ and
two in $A_{3-i}$, or has two vertices in $B_i$ and one in
$A_{3-i}$.
It follows that
\[
\tau_B(G)
=
\sum_{i=1}^2\sum_{u\in B_i}
\left(e(G[Y_u])+\frac12e_G(X_u,Y_u)\right)
+t(G[B]).
\]
Furthermore,
\[
\sum_{i=1}^2\sum_{u\in B_i}|X_u|
=
2e(G[B_1])+2e(G[B_2])
=O(b^2).
\]
Applying \eqref{eq:crude-triangles} to $G[B]$ gives
$t(G[B])=O_k(b^2)$, completing the proof.
\end{proof}

Define the missing-edge quantity
\[
D_B=
\sum_{i=1}^2\sum_{u\in B_i}
\bigl(a_{3-i}-d_{A_{3-i}}(u)\bigr).
\]
For $u\in E_i$, the summand is at least
$\eta m-|B_{3-i}|$.
Since $b=o(m)$, for sufficiently large $m$,
\begin{equation}\label{eq:DB-lower}
D_B
\ge
\frac{\eta}{2}m(|E_1|+|E_2|)
\ge \gamma_kmb,
\qquad
\gamma_k=\frac{\eta}{2K_k}>0.
\end{equation}
Consequently, Claim~\ref{clm:exceptional-triangles} yields
\begin{equation}\label{eq:exceptional-loss}
\tau_B(G)
\le
\frac r2(b_1a_2+b_2a_1)
-\frac r2\gamma_kmb
+O_k(b^2).
\end{equation}

Assume $
p=\Delta(G[A_1]),$ $s=\Delta(G[A_2]).$
By \eqref{eq:clean-properties},
$
0\le p,s\le r.
$


\begin{claim}\label{clm:degree-budget}
$p+s\le r.$
\end{claim}

\begin{proof}
Write $e_i=e(G[A_i])$.
We first prove that, if $p\ge 1$, then
\begin{equation}\label{eq:opposite-edge-bound}
e_2
\le
\frac{r-p}{2}a_2
+r(p+1)\eta m
+rp(r-p+2).
\end{equation}

Choose $u\in A_1$ with $d_{A_1}(u)=p$, and write
\[
N_G(u)\cap A_1=\{x_1,\ldots,x_p\}.
\]
Let
\[
Y=A_2\cap N_G(u)\cap\bigcap_{j=1}^pN_G(x_j).
\]
Since $u,x_1,\ldots,x_p\in T_1$,
\begin{equation}\label{eq:Y-missing}
|A_2\setminus Y|\le (p+1)\eta m.
\end{equation}

Let $\ell=r-p+2$.
Suppose that $G[Y]$ contains $p+1$ pairwise vertex-disjoint
copies of $P_\ell$.
Using $x_1,\ldots,x_p$ as successive connecting vertices,
we can concatenate these paths inside $G[N_G(u)]$.
The resulting path has
\[
(p+1)\ell+p=(p+1)(r-p+3)-1
\]
vertices. Moreover,
\begin{equation}\label{eq:concatenated-length}
(p+1)(r-p+3)-1-(2r+3)
=
(p-1)(r-p+1)\ge 0.
\end{equation}
Since $k\le 2r+3$, this contradicts
\eqref{eq:neighborhood-characterization}.

Choose a maximal family of pairwise vertex-disjoint
$P_\ell$'s in $G[Y]$, and let $Z$ be the union of their
vertex sets.
There are at most $p$ such paths, so
\[
|Z|\le p\ell.
\]
By maximality, $G[Y\setminus Z]$ is $P_\ell$-free.
Using Theorem \ref{thm8121},
$\Delta(G[A_2])\le r$, and \eqref{eq:Y-missing}, we obtain
\[
\begin{aligned}
e_2
&\le
e(G[Y\setminus Z])+r|Z|+r|A_2\setminus Y|\\
&\le
\frac{\ell-2}{2}|Y\setminus Z|
+rp\ell+r(p+1)\eta m\\
&\le
\frac{r-p}{2}a_2
+rp(r-p+2)+r(p+1)\eta m.
\end{aligned}
\]
This proves \eqref{eq:opposite-edge-bound}.

By symmetry, if $s\ge 1$, then
\begin{equation}\label{eq:symmetric-edge-bound}
e_1
\le
\frac{r-s}{2}a_1
+r(s+1)\eta m
+rs(r-s+2).
\end{equation}

Now suppose, for a contradiction, that $p+s\ge r+1$.
Since $p,s\le r$, both $p$ and $s$ are positive, so both
estimates apply.
The number of triangles entirely inside either $A_i$
is $O_k(m)$.
Every other triangle in $A_1\cup A_2$ consists of an
internal edge on one side and a vertex on the other.
Therefore
\[
\begin{aligned}
t(G[A_1\cup A_2])
&\le a_2e_1+a_1e_2+O_k(m)\\
&\le
\frac{2r-p-s}{2}a_1a_2
+r\eta m\bigl((s+1)a_2+(p+1)a_1\bigr)
+O_k(m)\\
&\le
\frac{r-1}{8}m^2+r(r+1)\eta m^2+O_k(m).
\end{aligned}
\]
By Claim~\ref{clm:exceptional-triangles} and $b=o(m)$,
\[
\tau_B(G)=O_k(bm+b^2)=o(m^2).
\]
Thus
\[
t(G)
\le
\left(\frac{r-1}{8}+r(r+1)\eta\right)m^2+o(m^2).
\]
By \eqref{eq:eta-conditions},
\[
t(G)
\le
\left(\frac r8-\frac3{32}\right)m^2+o(m^2).
\]
This contradicts
\[
t(G)\ge \gk(m)\ge \hk(m)=\frac r8m^2+O_k(m).
\]
Hence $p+s\le r$.
\end{proof}

Let
$M=a_1+a_2=m-b.$
The degree bounds give
\[
e(G[A_1])\le \frac p2a_1,\qquad
e(G[A_2])\le \frac s2a_2,
\]
and
\[
t(G[A_1])\le \frac{p(p-1)}6a_1,\qquad
t(G[A_2])\le \frac{s(s-1)}6a_2.
\]
Consequently,
\begin{equation}\label{eq:core-basic}
t(G[A_1\cup A_2])
\le
\frac{p+s}{2}a_1a_2
+\frac{p(p-1)}6a_1+\frac{s(s-1)}6a_2.
\end{equation}

\begin{claim}\label{lem:core-envelope}
\begin{equation}\label{eq:core-envelope}
t(G[A_1\cup A_2])
\le
\frac r8M^2+\frac{r(r-1)}{12} M+O_k(1).
\end{equation}
If $p+s=r$, then the stronger estimate
\begin{equation}\label{eq:core-mixed-loss}
t(G[A_1\cup A_2])
\le
\frac r8M^2+
\left(\frac{r(r-1)}{12}-\frac{ps}{6}\right)M+O_k(1)
\end{equation}
holds.
\end{claim}

\begin{proof}
If $p+s\le r-1$, \eqref{eq:core-basic} yields
\[
t(G[A_1\cup A_2])
\le \frac{r-1}{8}M^2+O_k(M).
\]
This implies \eqref{eq:core-envelope} for large $M$,
and a larger constant handles the bounded values of $M$.

Suppose that $p+s=r$, and write
\[
a_1=\frac M2+x,\qquad a_2=\frac M2-x.
\]
Then the right-hand side of \eqref{eq:core-basic} equals
\[
\frac r8M^2+
\left(\frac{r(r-1)}{12}-\frac{ps}{6}\right)M
-\frac r2x^2
+\frac{p(p-1)-s(s-1)}6x.
\]
The last two terms are bounded above by a constant
depending only on $k$.
This proves \eqref{eq:core-mixed-loss}, and
\eqref{eq:core-envelope} follows.
\end{proof}

\begin{claim}
$B_1=B_2=\varnothing$.
\end{claim}
\begin{proof}
  Combine 
\eqref{eq:exceptional-loss} and
Claim~\ref{lem:core-envelope}:
\[
\begin{aligned}
t(G)
\le{}&
\frac r8(m-b)^2+\frac{r(r-1)}{12}(m-b)
+\frac r2(b_1a_2+b_2a_1)\\
&-\frac r2\gamma_kmb+O_k(b^2)+O_k(1).
\end{aligned}
\]
By \eqref{eq:A-balanced},
\[
b_1a_2+b_2a_1=\frac{bm}{2}+o(mb).
\]
Hence
\begin{equation}\label{eq:elimination-bound}
t(G)
\le
\frac r8m^2+\frac{r(r-1)}{12} m
-\frac r2\gamma_kmb
+o(mb)+O_k(b^2)+O_k(1).
\end{equation}

If $b\ge 1$ along a sequence with $m\to\infty$, then
$b=o(m)$ implies $b^2=o(mb)$.
The negative term in \eqref{eq:elimination-bound}
therefore dominates the error terms and the bounded
difference between $\hk(m)$ and
$rm^2/8+\frac{r(r-1)}{12} m$.
This contradicts $t(G)\ge \gk(m)\ge \hk(m)$.
Consequently, for all sufficiently large $m$,
\begin{equation}\label{eq:B-empty}
B_1=B_2=\varnothing.
\end{equation}
\end{proof}

If $p,s\ge 1$, then \eqref{eq:core-mixed-loss} gives
\[
t(G)
\le
\frac r8m^2+
\left(\frac{r(r-1)}{12}-\frac{ps}{6}\right)m+O_k(1),
\]
again contrary to \eqref{eq:lower-bound}.
After interchanging the parts, we may therefore assume
\[
p=r,\quad s=0.
\]
Thus now we have
$
V(G)=A_1\dotcup A_2,
$
$A_2$ is independent, and
$H=G[A_1]$ satisfies $\Delta(H)= r$ and is $P_k$-free.
It follows that
\begin{equation}\label{eq:final-g-bound}
t(G)
\le |A_2|e(H)+t(H)
\le \gk(m).
\end{equation}
Combining this with \eqref{eq:deletion-surplus}, we obtain
\[
\gk(m)+n-m\le t(G)\le \gk(m).
\]
Hence $m=n$, and
\[
t(G_n)\le \gk(n).
\]
We conclude that
\begin{equation}\label{eq:ex-equals-g}
\ex(n, K_3, \widehat{P_k})=\gk(n)
\end{equation}
for all sufficiently large $n$.
\begin{lemma}\label{lem:regular}
For fixed $k\ge 6$ and sufficiently large $n$, every graph
$H$ attaining the maximum in
\eqref{eq:g-definition-intro} is $r$-regular.
\end{lemma}

\begin{proof}
Put $a=v(H)$, and define
\[
D_e=\frac{ra}{2}-e(H),\qquad
D_t=\frac{r(r-1)}6a-t(H).
\]
Both quantities are nonnegative because $\Delta(H)\le r$.
Using \eqref{eq:complete-square}, we obtain the exact identity
\[
(n-a)e(H)+t(H)
=
U(n)-\frac r2(a-x_n)^2-(n-a)D_e-D_t.
\]
Since $H$ attains $\gk(n)$ and
$\gk(n)\ge \hk(n)\ge U(n)-rq^2/8$ for large $n$,
\begin{equation}\label{eq:regularity-budget}
\frac r2(a-x_n)^2+(n-a)D_e+D_t
\le \frac{rq^2}{8}.
\end{equation}
In particular,
\[
a=\frac n2+O_k(1).
\]
If $D_e>0$, then $2D_e$ is a positive integer, so
$D_e\ge 1/2$.
The term $(n-a)D_e$ would then tend to infinity, contradicting
\eqref{eq:regularity-budget}.
Therefore $D_e=0$.
Since every degree is at most $r$ and the average degree is
$r$, all degrees are equal to $r$.
\end{proof}




Lemma~\ref{lem:join} proves the converse assertion in
Theorem~\ref{thm:main-structure}, so the theorem is now
proved in full.

\section{Proof of Theorem \ref{main2}}


For an $r$-regular graph $H$, recall
\[
D(H)=\frac{r(r-1)}6v(H)-t(H).
\]
Since
\[
3t(H)=\sum_{v\in V(H)}e(H[N_H(v)])
\le v(H)\binom r2,
\]
we have $D(H)\ge 0$.
Moreover,
\begin{equation}\label{eq:zero-deficit}
D(H)=0
\quad\Longleftrightarrow\quad
H\text{ is a disjoint union of copies of }K_q.
\end{equation}
Indeed, equality forces the neighborhood of every vertex
to be a clique.
Since $H$ is $r$-regular, the closed neighborhood of each
vertex then forms a $K_q$ component.
The converse is immediate.
Recall
$$
x_n=\frac{n}{2}+\frac{r-1}{6}, \quad U(n)=\frac{r}{8} n^2+\frac{r(r-1)}{12} n+\frac{r(r-1)^2}{72} .
$$
For $a=v(H)$, we have
\begin{equation}\label{eq:join-deficit}
t(H+\Ik{n-a})
=
U(n)-
\left[
\frac r2(a-x_n)^2+D(H)
\right].
\end{equation}

For an $r$-regular graph $H$, define its triangle deficit by
$$
D(H)=\frac{r(r-1)}{6} v(H)-t(H) .
$$
For $0 \leq s<q$, define
$$
d_s=\min \left\{D(F): \begin{array}{l}
F \text { is } r \text {-regular and } P_k \text {-free, } \\
F \text { has no connected component isomorphic to } K_q, \\
v(F) \leq \left\lfloor\frac{3 r q^2}{8}\right\rfloor , \quad v(F) \equiv s \quad(\bmod q)
\end{array}\right\} .
$$
The minimum of an empty set is $+\infty$. The empty graph is allowed and is regarded as $r$-regular vacuously; in particular, $d_0=0$.

\begin{claim}\label{thm:T-formula}
For fixed $k\ge 6$ and all sufficiently large $n$,
\begin{equation}\label{eq:T-formula}
ex(n, K_3, \widehat{P_k})
=
U(n)-
\min_{\substack{0\le s<q\\d_s<\infty}}
\left\{
\frac r2\dist(x_n,s+q\Z)^2+d_s
\right\}.
\end{equation}
\end{claim}

\begin{proof}

Let $H+\Ik{n-a}$ attain $\gk(n)$.
Then by (\ref{eq:h-envelope}), (\ref{eq:join-deficit}) and $\hk(n)\leq\gk(n)$, we have
\begin{equation}\label{eq:defect-budget}
\frac r2(a-x_n)^2+D(H)\le \frac{rq^2}{8}.
\end{equation}
Delete all $K_q$ components of $H$, and let $F$ be the
remaining graph. It is $r$-regular and $P_k$-free, has no
$K_q$ component, and
\[
D(F)=D(H).
\]

For every vertex $v\in V(F)$,
\[
\binom r2-e(F[N_F(v)])\ge 1.
\]
Otherwise $N_F(v)$ would be a clique, and regularity would
make $N_F[v]$ a $K_q$ component.
Summing gives
\[
v(F)
\le
\sum_{v\in V(F)}
\left(\binom r2-e(F[N_F(v)])\right)
=
3D(F).
\]
By \eqref{eq:defect-budget},
\[
v(F)\le 3D(H)\le \frac{3rq^2}{8}.
\]

Let $a\equiv s\pmod q$, then
\[
v(F)\equiv s\pmod q.
\]
Thus
\[
D(H)=D(F)\ge d_s,
\qquad
|a-x_n|\ge \dist(x_n,s+q\Z).
\]
Using \eqref{eq:join-deficit}, we obtain
\[
ex(n, K_3, \widehat{P_k})
\le
U(n)-
\min_{\substack{0\le s<q\\d_s<\infty}}
\left\{
\frac r2\dist(x_n,s+q\Z)^2+d_s
\right\}.
\]

For the reverse inequality, fix $s$ with $d_s<\infty$,
and choose a graph $F$ attaining $d_s$.
Let $a\in s+q\Z$ be nearest to $x_n$.
Since $v(F)\le \left\lfloor\frac{3 r q^2}{8}\right\rfloor $ depends only on $k$, for sufficiently
large $n$ we have
\[
v(F)\le a\le n.
\]
Define
\[
H=\frac{a-v(F)}qK_q\dotcup F.
\]
Then $H$ is $r$-regular and $P_k$-free, and $D(H)=d_s$.
Equation~\eqref{eq:join-deficit} gives
\[
t(H+\Ik{n-a})
=
U(n)-
\left\{
\frac r2\dist(x_n,s+q\Z)^2+d_s
\right\}.
\]
Choose a residue minimizing the bracketed expression.
This proves the reverse inequality and hence
\eqref{eq:T-formula}.
\end{proof}


Let
\[
c_k(n)
=
\frac{r(r-1)^2}{72}
-
\min_{\substack{0\le s<q\\d_s<\infty}}
\left\{
\frac r2\dist(x_n,s+q\Z)^2+d_s
\right\}.
\]
The minimum is finite because $d_0=0$.
It is bounded below by $0$ and above by $rq^2/8$.
Thus $c_k(n)$ is bounded.
Since
\[
x_{n+2q}=x_n+q,
\]
the distance from $x_n$ to each lattice $s+q\Z$ is
unchanged when $n$ is replaced by $n+2q$.
Therefore $c_k(n+2q)=c_k(n)$.
From
Claim~\ref{thm:T-formula}, we have \[
ex(n, K_3, \widehat{P_k})
=
\frac r8n^2+\frac{r(r-1)}{12}n+c_k(n).
\]
Theorem~\ref{main2} now follows immediately.

\vspace{5pt}
\noindent{\bf Declaration of AI assistance}:
During the preparation of this work, the authors used OpenAI ChatGPT to assist with the development of some
mathematical proofs. The authors reviewed the AI-assisted material and take full responsibility for the correctness
of the proofs and the content of the article.

\vspace{5pt}
\noindent{\bf Acknowledgements}:
This work was supported by the National Key Research and Development Program of China (2023YFA1010203), and the Innovation Program for Quantum Science and Technology (2021ZD0302902).

\bibliographystyle{plain}  
\bibliography{ref}  

\end{document}